\documentclass[sts]{imsart}

\RequirePackage[OT1]{fontenc}
\usepackage{amsthm,amsmath,natbib}
\RequirePackage[colorlinks,citecolor=blue,urlcolor=blue]{hyperref}
\usepackage{graphicx}

\startlocaldefs
\numberwithin{equation}{section}
\theoremstyle{plain}

\newtheorem{Result}{Result}[section]
\endlocaldefs

\begin{document}

\begin{frontmatter}
\title{Are profile likelihoods likelihoods? No, but sometimes they can be.}
\runtitle{Profile likelihood and equivalent priors}

\begin{aug}
\author{\fnms{Alan} \snm{Huang}\ead[label=e1]{alan.huang@uq.edu.au}},
\and
\author{\fnms{Andy Sangil} \snm{Kim}\ead[label=e2]{andykim86@gmail.com}}

\runauthor{A. Huang \& A. S. Kim}

\affiliation{University of Queensland, and University of Technology Sydney}

\address{School of Mathematics and Physics, University of Queensland, St Lucia, QLD, Australia 4072 \printead{e1}.}

\address{School of Mathematical and Physical Sciences, University of Technology Sydney, Ultimo, NSW, Australia 2007 \printead{e2}.}
\end{aug}

\begin{abstract}
We offer our two cents to the ongoing discussion on whether profile likelihoods are ``true" likelihood functions, by showing that the profile likelihood function can in fact be identical to a marginal likelihood in the special case of normal models. Thus, profile likelihoods can be ``true" likelihoods insofar as marginal likelihoods are ``true" likelihoods. The prior distribution that achieves this equivalence turns out to be the Jeffreys prior. We suspect, however, that normal models are the only class of models for which such an equivalence between maximization and marginalization is exact.
\end{abstract}


\end{frontmatter}

\section{Introduction}
There has been some recent and rather lively debate as to whether the profile likelihood, obtained by maximizing out nuisance parameters in the full likelihood, can be considered a ``true" likelihood function in the remaining parameters, with arguments ranging from probabilistic, possibilistic and even philosophical perspectives \citep[e.g.,][]{Aitkin2005, Aitkin2010, Evans2015, Maclaren2018, Robert2018}. 
Here, the notion of a ``true" likelihood is a function that corresponds to some joint probability distribution on the data for each value of the model parameters\footnote{\cite{Maclaren2018} argues for a different notion of likelihood based on possibility rather than probability, with addition replaced by maximization. This approach is worth further consideration, however in this note we stick with the classical probability-based notion of likelihood.}.

The consensus from the statistical literature
seems to be ``no", in general. \cite{Aitkin2005} states rather unequivocally that {\it ``the profile likelihood is not a likelihood, but a likelihood maximized over nuisance parameters given the values of the parameters of interest."} 
In other words,
the maximization operator
does not generally take probability distributions to probability distributions, but merely to a ``slice" in a probability distribution (hence the ``profile" moniker).

Of course, from a frequentist point of view the profile likelihood can still exhibit likelihood-type statistical properties, regardless of whether or not it corresponds to a true likelihood. These properties include consistency, asymptotic normality and asymptotic efficiency of its maximizer, with the profile likelihood ratio test even exhibiting Wilks' phenomenon under some general conditions \citep{MVDV2000}.

From a Bayesian point of view, nuisance parameters are usually dealt with via marginalization instead of maximization. In contrast to the profile likelihood, there is little debate as to whether the marginal likelihood corresponds to a true likelihood, as {\it ``integration over variables takes probability distributions to probability distributions"} \citep{Maclaren2018}. Indeed, the elementary concept of a marginal probability is constructed precisely by integrating joint probabilities over a subset of variables.

While maximization and marginalization are two seemingly disparate operators, it turns out that in the special case of normal models, the profile likelihood for the mean parameter(s) is precisely equivalent to the marginal likelihood obtained by integrating over Jeffreys prior on the nuisance variance parameters. In this case, profile likelihood can be considered a true likelihood insofar as a marginal likelihood is a true likelihood. This equivalence is exact for normal models, and we speculate that, like other results from likelihood theory, it may only be asymptotically true for other exponential families. 


\section{Profile likelihood, marginal likelihood and the equivalent prior for normal models}
Let $y = (y_1, y_2, \ldots, y_n)^\top$ be a random sample from a normal distribution with mean $\mu$ and variance $\sigma^2$. The likelihood function (using Bayesian notation) is given by
$$
p(y|\mu, \sigma^2) \propto (\sigma^2)^{-n/2} \exp\left\{- \frac{1}{2\sigma^2} \sum_{i=1}^n (y_i-\mu)^2 \right\} \ .
$$
The maximum likelihood estimator of $\sigma^2$ for each given $\mu$ is
$$
\hat \sigma^2(\mu) = \frac{1}{n} \sum_{i=1}^n (y_i -\mu)^2 \ ,
$$
so that the profile likelihood for $\mu$ is
$$
\sup_{\sigma^2 } p(y|\mu, \sigma^2) \, = \, p(y|\mu, \hat \sigma^2(\mu)) \, \propto \, \left[\sum_{i=1}^n (y_i -\mu)^2 \right]^{-n/2} \ .
$$
It is not immediately clear that this function corresponds to a valid probability distribution in $y$ for each $\mu$. This kind of ambiguity is precisely what has fuelled the debate over whether profile likelihoods can be considered true likelihoods.

On the other hand, consider a Jeffreys prior $p(\sigma^2) \propto 1/\sigma^2$ on the variance $\sigma^2$ where the mean $\mu$ is treated as given. Integrating out $\sigma^2$ leads to the marginal likelihood as
\begin{eqnarray*}
p(y|\mu) &=& \int_{\sigma^2} p(y|\mu, \sigma^2) \, p(\sigma^2) \, d(\sigma^2)  
\\
 &\propto &  \int_{\sigma^2} (\sigma^2)^{-(n+2)/2} \exp\left\{ - \frac{1}{2\sigma^2} \sum_{i=1}^n (y_i - \mu)^2 \right\} d(\sigma^2)  
\\
& \propto& \left[ \sum_{i=1}^n (y_i - \mu)\right]^{-n/2} \ ,
\end{eqnarray*}
by noticing that the integrand is the kernel of a Inverse-Gamma distribution with shape parameter $n/2$ and scale parameter $\sum_{i=1}^n (y_i - \mu)^2/2$. We see that the marginal likelihood coincides exactly with the profile likelihood, that is,
$$
\sup_{\sigma^2} p(y|\mu,\sigma^2) \equiv \int_{\sigma^2} p(y|\mu, \sigma^2) \, p(\sigma^2) \, d(\sigma^2)  
$$
for Jeffreys prior $p(\sigma^2) \propto 1/\sigma^2$ on the variance $\sigma^2$. 

\begin{figure}
\includegraphics[width = 0.6\textwidth]{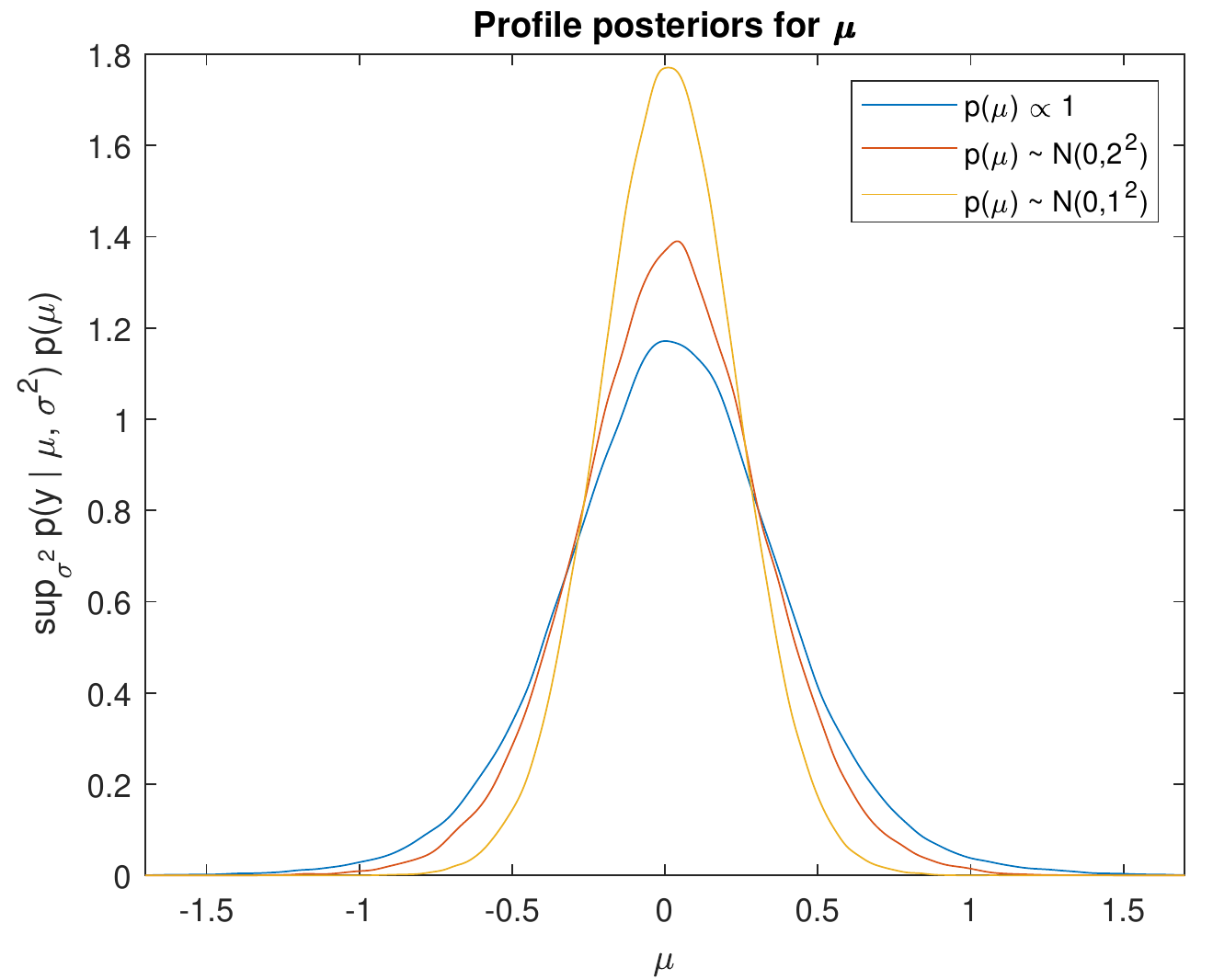}
\caption{Profile posterior distributions for $\mu$ obtained by multiplying the profile likelihood by various prior distributions $p(\mu)$ given a random sample of $n=10$ standard normal observations.}
\label{fig:profilepost}
\end{figure}

A practical consequence is that the profile likelihood can be used to construct valid posterior distributions for Bayesian inferences. Given a prior $p(\mu)$ on $\mu$, the `profile posterior" is
$$
 \sup_{\sigma^2} p(y | \mu, \sigma^2) \, p(\mu) \propto  \int_{\sigma^2}  p(y | \mu, \sigma^2) \, p(\mu) \,  p(\sigma^2) d(\sigma^2) \ ,
$$
precisely the same as the marginal posterior obtained from integrating over
Jeffreys prior on $\sigma^2$. For example, profile posterior distributions resulting from priors $p(\mu) \sim N(0,1^2), p(\mu) \sim N(0,2^2),$ and the improper prior $p(\mu) \propto 1$ can be obtained via Gibbs sampling, following Chapter 8.2.1 of \citet{KC2014}, say. Indeed, following the example in that book, Figure \ref{fig:profilepost} displays the corresponding profile posterior distributions for $\mu$ given a random sample of 10 observations from a standard normal distribution.

The Jeffreys prior can therefore be thought of as an ``equivalent prior" that makes marginalizing the likelihood equivalent to maximizing the likelihood. Analogous results for the multivariate and regression cases are also elementary to show.
\begin{Result}
Let $y_1,y_2,\ldots,y_n \stackrel{\rm{iid}}{\sim} N_d(\mu, \Sigma)$, where $\mu \in \mathbf{R}^d$ is a $d$-vector of mean parameters of interest and $\Sigma$ is a $d \times d$ nuisance variance matrix. Then the profile likelihood for $\mu$ is 
equivalent to the marginal likelihood for $\mu$
for Jeffreys prior $p(\Sigma) \propto |\Sigma|^{-(d+1)/2}$ on $\Sigma$.
\label{pr:iidvector}
\end{Result}

\begin{Result}
Let $y_i|x_i \stackrel{\rm ind}{\sim} N(x_i^\top \beta, \sigma^2)$, $i=1,2,\ldots,n$, where each $x_i  \in \mathbf{R}^q$ is a vector of covariates, $\beta$ is an associated vector of mean parameters of interest and $\sigma^2$ is a nuisance variance parameter. Then the profile likelihood for $\beta$ is 
equivalent to the marginal likelihood for $\beta$
for Jeffreys prior $p(\sigma^2) \propto 1/\sigma^2$ on $\sigma^2$.
\end{Result}

\section{Discussion}
Our contribution to the ongoing debate over the nature of the profile likelihood is to provide a simple (counter-)example in which the profile likelihood is 
identical to a marginal likelihood. We find it rather remarkable and somewhat counter-intuitive that marginalization can be made equivalent to maximization via a particular choice of prior on the nuisance parameters. That this equivalent prior happens to be the well-known Jeffreys prior is also an interesting coincidence, but perhaps not completely unexpected as both the profile likelihood and Jeffreys prior are constructed to be ``non-informative" in some frequentist or Bayesian sense, respectively. Of course, whether the improper Jeffreys prior constitutes a ``true" prior 
that can be integrated over is another debate, perhaps for another day.

In keeping with \cite{Aitkin2005}, we suspect that normal models are the only class of models for which this equivalence is exact. However, we also speculate that a generalization may hold asymptotically for other exponential families. Heuristically speaking, the score equations for exponential families, whilst typically not solvable in closed-form, can be linearized in its parameters, with leading term proportional to the Hessian of the likelihood, the inverse of which forms the basis of Jeffreys prior. It is also well-known that exponential families for data induce exponential families in the model parameters (which is why exponential families always have conjugate priors). These two ingredients combine to give us hope that the (linearized) profile likelihood might pop up as the normalizing constant when integrating out an exponential family likelihood over the Jeffreys prior, just as it did in the normal case. This is a lead worth exploring further.

\section*{Acknowledgements}
We thank Dr Yao-ban Chan (Melbourne) for comments that improved this note.


\begin{thebibliography}{}

\bibitem[\protect\citeauthoryear{Aitkin}{2005}]{Aitkin2005}
\textsc{Aitkin, M.} (2005). Profile Likelihood. In \textit{Encyclopedia of Biostatistics}, John Wiley \& Sons.

\bibitem[\protect\citeauthoryear{Aitkin}{2010}]{Aitkin2010}
\textsc{Aitkin, M.} (2010). \textit{Statistical Inference: An Integrated Bayesian/Likelihood Approach}, Chapman \& Hall/CRC Monographs on Statistics \& Applied Probability. CRC Press.

\bibitem[\protect\citeauthoryear{Evans}{2015}]{Evans2015}
\textsc{Evans, M.} (2015). \textit{Measuring Statistical Evidence Using Relative Belief}, Chapman \& Hall/CRC Monographs on Statistics \& Applied Probability. CRC Press.

\bibitem[\protect\citeauthoryear{Kroese \& Chan}{2014}]{KC2014}
\textsc{Kroese, D.P. \& Chan, J.C.C.} (2014). \textit{Statistical Modeling and Computation}, Springer, NY.

\bibitem[\protect\citeauthoryear{Maclaren}{2018}]{Maclaren2018}
\textsc{Maclaren, O.J.} (2018).
Is profile likelihood a true likelihood? An argument in favor. \textit{arXiv arxiv.org/abs/1801.04369}

\bibitem[\protect\citeauthoryear{Murphy \& van der Vaart}{2000}]{MVDV2000}
\textsc{Murphy, S.A \& van der Vaart, A. W.} (2000). On Profile Likelihood. \textit{Journal of the American Statistical Association}, \textbf{95}, 449--465.

\bibitem[\protect\citeauthoryear{Robert}{2018}]{Robert2018}
\textsc{Robert, C.P.} (2018). \textit{xianblog.wordpress.com/2018/03/27/are-profile-likelihoods-likelihoods/}

\end{thebibliography}
\end{document}